\documentclass[11pt]{amsart}

\usepackage{amsmath}
\usepackage{amssymb}
\usepackage{graphicx}
\theoremstyle{definition}

\numberwithin{equation}{section}
\title[Persistence of pulses for certain reaction-diffusion equations]
{Persistence of pulses for certain reaction-diffusion equations in dimensions
two and three}

\author[V. Vougalter]{Vitali Vougalter}

\address[V. Vougalter]{Department of Mathematics,  University of Toronto,
Toronto, Ontario, M5S 2E4, Canada}
\email{{\tt vitali@math.toronto.edu}}



\begin{document}

\begin{abstract}
We address the persistence under a perturbation of
stationary pulse solutions of some reaction-diffusion type equations in
dimensions $d=2,3$ and evaluate the asymptotic approximations of such pulses
to the leading order in the parameter of the perturbation.
\end{abstract}

\maketitle

%%%%%%%%%%%%%%%%%%%%%%%%%%%%%%%%%%%%%%%%%%%%%%%%%%%%%%%%%%%%%%%%%%%%%%%%%%%%%%%%

\setcounter{equation}{0}

\section{Introduction}

\noindent
Consider the equation,
\begin{equation}
\label{rd}
\frac{\partial u}{\partial t}=\Delta u+au^{2}
(1-u)-\sigma(x)u, \quad x\in {\mathbb R}^{d}, \quad d=2,3, \quad t\geq 0.
\end{equation}
Here $a>0$ is a constant. Problems of this kind are important for the study of 
the population dynamics (see e.g. ~\cite{V14}).
In such context, $u(x,t)$ stands for the density of a population at time $t$
in location $x$, the diffusion
term describes its migration, the second term in the right side is the
reproduction rate, and the last term is the mortality rate. In the situation of
the sexual reproduction, the reproduction rate is proportional to the second
power of the population density and to the available resources $(1-u)$, given
by the difference of the rate of production of resources and the rate of their
consumption. Stationary solutions of (\ref{rd}) solve
\begin{equation}
\label{st}
\Delta u+au^{2}
(1-u)-\sigma(x)u=0, \quad x\in {\mathbb R}^{d}, \quad d=2,3.
\end{equation}
Stationary solutions of reaction-diffusion type problems decaying at 
infinities are called
{\it pulses}. Their existence along with the stability and other related issues
have been studied actively for both local and nonlocal models in recent years,
for example, in ~\cite{AM07}, ~\cite{BLP81}, ~\cite{CV21}, ~\cite{EVV17},
~\cite{IW00}, ~\cite{KP97}, ~\cite{MV23}, 
~\cite{TBMV11}, ~\cite{VV141}, ~\cite{WW07}, ~\cite{Z04}.
In the present work, we only focus on the persistence of pulses and
their asymptotic approximations to the leading order of the small parameter
when a perturbation is applied to problem ~(\ref{rd}). Similar ideas in the
context of standing solitary waves of the nonlinear Schr\"odinger equation,
when a perturbation is applied to either a scalar potential contained in it or
to the nonlinear term, were exploited in ~\cite{V10}.

The unperturbed stationary equation in our case is 
\begin{equation}
\label{st0}
\Delta w_{0}+aw_{0}^{2}
(1-w_{0})-\sigma_{0}(x)w_{0}=0, \quad x\in {\mathbb R}^{d}, \quad d=2,3.
\end{equation}
In the work we consider the space $H^{2}({\mathbb R}^{d})$ equipped
with the norm
\begin{equation}
\label{h2n}
\|u\|_{H^{2}({\mathbb R}^{d})}^{2}:=\|u\|_{L^{2}({\mathbb R}^{d})}^{2}+ \|
\Delta u\|_{L^{2}({\mathbb R}^{d})}^{2}, \quad d=2,3.
\end{equation}
By means of the standard Sobolev embedding, we have
\begin{equation}
\label{em}
\|u\|_{L^{\infty}({\mathbb R}^{d})}\leq c_{e}\|u\|_{H^{2}({\mathbb R}^{d})}, \quad d=2,3,
\end{equation}
where $c_{e}>0$ is the constant of the embedding.
Let us first state the following conditions on the parameters of our unperturbed
equation along
with its solution that we are going to consider.

\bigskip

\noindent
{\bf Assumption 1.1.} {\it Let $d=2,3$, the constant $a>0$ and the function}
$$
\sigma_{0}(x): {\mathbb R}^{d}\to {\mathbb R}, \ 
\sigma_{0}(x)\in C^{\infty}({\mathbb R}^{d}), \ \hbox{lim}_{|x|\to \infty}\sigma_{0}
(x)=\delta>0.
$$
{\it We also assume that problem (\ref{st0}) has a pulse solution
$w_{0}(x)>0, \ x\in{\mathbb R}^{d}$, which satisfies 
\begin{equation}
\label{w0}  
w_{0}(x)\in C^{\infty}({\mathbb R}^{d}), \ w_{0}(x)\leq Ce^{-\alpha |x|}, \
x\in{\mathbb R}^{d},
\end{equation}
where $C,\alpha>0$ are the constants.}

\bigskip

It can be easily deduced from equation (\ref{st0}) via Assumption 1.1 that
$w_{0}(x)\in H^{2}({\mathbb R}^{d})$ because it decays
exponentially as $|x|\to \infty$. Here we
used the conditions on $\sigma_{0}(x)$.

When a perturbation is applied to stationary problem (\ref{st0}), we
have
\begin{equation}
\label{stp}
\Delta w+aw^{2}
(1-w)-[\sigma_{0}(x)+\varepsilon \sigma_{1}(x)]w=0, \quad x\in {\mathbb R}^{d},
\quad d=2,3.
\end{equation} 

\bigskip

\noindent
{\bf Assumption 1.2.} {\it Let the parameter $\varepsilon\geq 0$ and the
nontrivial function}
\begin{equation}
\label{ass2}
\sigma_{1}(x): {\mathbb R}^{d}\to {\mathbb R}, \  
\sigma_{1}(x)\in C^{\infty}({\mathbb R}^{d}), \ \hbox{lim}_{|x|\to \infty}
\sigma_{1}(x)=0.
\end{equation}

\bigskip

Let us look for the solutions of problem (\ref{stp}) in the form
\begin{equation}
\label{0p}
w(x)=w_{0}(x)+w_{p}(x).
\end{equation}
By virtue of (\ref{stp}) along with (\ref{st0}), we arrive at
\begin{equation}
\label{pert}
L_{0}w_{p}(x)=a(1-3w_{0}(x)){w_{p}^{2}(x)}-a{w_{p}^{3}(x)}-\varepsilon \sigma_{1}(x)
(w_{0}(x)+w_{p}(x))
\end{equation}
with
\begin{equation}
\label{l0}
L_{0}=-\Delta+a(3w_{0}^{2}(x)-2w_{0}(x))+\sigma_{0}(x):
H^{2}({\mathbb R}^{d})\to L^{2}({\mathbb R}^{d}), \ d=2,3.
\end{equation}
Under Assumption 1.1, it can be trivially checked that the essential spectrum
of $L_0$ is
\begin{equation}
\label{ess}
\sigma_{ess}(L_{0})=[\delta, +\infty).
\end{equation}
If $\sigma_{0}(x)$ had been a constant
function in ${\mathbb R}^{d}$, then the operator $L_{0}$ would have had zero
modes
$\displaystyle{\frac{\partial w_{0}}{\partial x_{k}}, \ 1\leq k\leq d}$. This
can be easily obtained by differentiating
both sides of problem (\ref{st0}). However, in the present work we assume
the function $\sigma_{0}(x)$ to be generic such that operator (\ref{l0})
would have a trivial kernel.

\bigskip

\noindent
{\bf Assumption 1.3.} {\it The kernel} $\hbox{ker}(L_{0})=\{0\}$.

\bigskip

By virtue of (\ref{ess}) along with Assumption 1.3, the operator
\begin{equation}
\label{l2h2}  
L_{0}^{-1}: L^{2}({\mathbb R}^{d})\to H^{2}({\mathbb R}^{d}), \quad d=2,3
\end{equation}
is bounded, so that its norm
\begin{equation}
\label{l0-}
\|L_{0}^{-1}\|<\infty.
\end{equation}
We introduce a closed ball in the Sobolev space $H^2(\mathbb{R}^{d})$ as
\begin{equation}
\label{b}
B_{\rho}:=\{u\in H^{2}({\mathbb R}^{d}) \ | \ \|u\|_{H^{2}({\mathbb R}^{d})}\leq \rho \},
\quad \rho>0, \quad d=2,3.
\end{equation}
Let us seek the solutions of equation (\ref{pert}) as fixed points of the
auxiliary nonlinear problem
\begin{equation}
\label{aux}
L_{0}u(x)=a(1-3w_{0}(x)){v^{2}(x)}-a{v^{3}(x)}-\varepsilon \sigma_{1}(x)
(w_{0}(x)+v(x)).
\end{equation}
For a given function $v(x)$ this is an equation with respect to $u(x)$.
Note that the   
similar ideas for problems containing non-Fredholm operators in their
left sides have been exploited in ~\cite{EV24} and ~\cite{VV24}.
Existence of solutions of certain quadratic integral equations in dimensions
two and three was discussed in ~\cite{V24}.
Let us
introduce the operator $t$, such that $u=tv$, where $u$ is a solution of
equation (\ref{aux}).
Our main statement is as follows.

\bigskip

\noindent
{\bf Theorem 1.4.} {\it Let Assumptions 1.1, 1.2 and 1.3 be valid. Then
problem (\ref{aux}) defines the map $t: B_{\rho}\to B_{\rho}$, which is a strict
contraction for all $0<\rho<\rho^{*}$ and $0<\varepsilon<\varepsilon^{*}$ for
certain $\rho^{*}>0$ and $\varepsilon^{*}>0$. The unique fixed point $w_{p}(x)$ of
this map $t$ is the only solution of equation (\ref{pert}) in $B_{\rho}$, so
that
\begin{equation}
\label{ass}
w_{p}(x)=-\varepsilon L_{0}^{-1}[\sigma_{1}(x)w_{0}(x)]+O({\varepsilon}^{2}).
\end{equation}
The resulting solution of problem (\ref{stp}) \ $w(x)\in
H^{2}({\mathbb R}^{d})$ is given by formula (\ref{0p}).}

\bigskip

We note that $O({\varepsilon}^{2})$ in the right side of (\ref{ass})
stands for the terms of the order ${\varepsilon}^{2}$ and higher in the
$H^{2}({\mathbb R}^{d})$ norm. The proof of Theorem 1.4 is presented in the 
following section.

\bigskip 

%%%%%%%%%%%%%%%%%%%%%%%%%%%%%%%%%%%%%%%%%%%

\setcounter{equation}{0}

\section{Proof of Theorem 1.4}

\noindent
{\it Proof.} Let us first demonstrate the uniqueness of
solutions of equation (\ref{aux}).
Suppose that there exists  a $v\in B_{\rho}$, so that ~(\ref{aux}) admits two
solutions 
$u_1$, $u_2\in B_{\rho}$.   Then the difference function
$\xi(x):=u_{1}(x)- u_{2}(x)\in H^{2}({\mathbb R}^{d})$ is a solution of the
homogeneous problem
\begin{equation}
\label{hom}
L_{0}\xi=0.
\end{equation}
By virtue of Assumption 1.3, equation (\ref{hom}) has only a trivial solution.
This gives us the uniqueness of solutions of ~(\ref{aux}).

Then, for an arbitrary $v(x)\in B_{\rho}$, we obtain the upper bound on the
right side of equation (\ref{aux}) in the absolute value as
$$
[a(1+3\|w_{0}\|_{L^{\infty}({\mathbb R}^{d})})\|v\|_{L^{\infty}({\mathbb R}^{d})}+
$$
\begin{equation}
\label{ub1}
+a\|v\|_{L^{\infty}({\mathbb R}^{d})}^{2}+\varepsilon \|\sigma_{1}\|_
{L^{\infty}({\mathbb R}^{d})}]
|v(x)|+\varepsilon \|\sigma_{1}\|_{L^{\infty}({\mathbb R}^{d})}|w_{0}(x)|.
\end{equation}
Clearly, $\sigma_{1}(x)\in L^{\infty}({\mathbb R}^{d})$ via Assumption 1.2. By
virtue of the Sobolev embedding (\ref{em}), quantity (\ref{ub1}) can be
estimated from above by
$$
[a(1+3c_{e}\|w_{0}\|_{H^{2}({\mathbb R}^{d})})c_{e}\|v\|_{H^{2}({\mathbb R}^{d})}+
$$
\begin{equation}
\label{ub2}
+ac_{e}^{2}\|v\|_{H^{2}({\mathbb R}^{d})}^{2}+\varepsilon \|\sigma_{1}\|_
{L^{\infty}({\mathbb R}^{d})}]
|v(x)|+\varepsilon \|\sigma_{1}\|_{L^{\infty}({\mathbb R}^{d})}|w_{0}(x)|.
\end{equation}
Since $v(x)\in B_{\rho}$, we have the upper bound for (\ref{ub2}) equal to
$$
[a(1+3c_{e}\|w_{0}\|_{H^{2}({\mathbb R}^{d})})c_{e}\rho+
$$
\begin{equation}
\label{ub3}
+ac_{e}^{2}\rho^{2}+\varepsilon \|\sigma_{1}\|_{L^{\infty}({\mathbb R}^{d})}]
|v(x)|+\varepsilon \|\sigma_{1}\|_{L^{\infty}({\mathbb R}^{d})}|w_{0}(x)|.
\end{equation}
Using (\ref{aux}), we easily derive that
$$
\|u\|_{H^{2}({\mathbb R}^{d})}\leq \|L_{0}^{-1}\|\{[a(1+3c_{e}\|w_{0}\|_
{H^{2}({\mathbb R}^{d})})c_{e}\rho+
$$
\begin{equation}
\label{ub4}
+ac_{e}^{2}\rho^{2}+\varepsilon \|\sigma_{1}\|_{L^{\infty}({\mathbb R}^{d})}]\rho
+\varepsilon \|\sigma_{1}\|_{L^{\infty}({\mathbb R}^{d})}\|w_{0}\|_{H^{2}({\mathbb R}^{d})}\}.
\end{equation}
Evidently, the estimate
$$
\|L_{0}^{-1}\|\{ac_{e}(1+3c_{e}\|w_{0}\|_{H^{2}({\mathbb R}^{d})})\rho^{2}+
$$
\begin{equation}
\label{ub5}
+ac_{e}^{2}\rho^{3}+\varepsilon \|\sigma_{1}\|_{L^{\infty}({\mathbb R}^{d})}(\rho+
\|w_{0}\|_{H^{2}({\mathbb R}^{d})})\}\leq \rho
\end{equation}
can be achieved for all $\rho>0$ and $\varepsilon>0$ sufficiently small.
Clearly, 
the upper bound on the values of $\varepsilon>0$ here will depend on $\rho$.
Therefore,
\begin{equation}
\label{ur}
\|u\|_{H^{2}({\mathbb R}^{d})}\leq \rho,
\end{equation}
so that  $u\in B_{\rho}$ as well. Thus, problem (\ref{aux}) defines the
map $t: B_{\rho}\to B_{\rho}$ for all $\rho>0$ and $\varepsilon>0$ 
small enough.

The goal is to show that this map is a strict contraction. Let us choose
arbitrarily $v_1$, $v_2\in B_{\rho}$. The reasoning above implies that
$u_1:=tv_1$, $u_2:=tv_2 \in B_{\rho}$ as well for all $\rho>0$ and
$\varepsilon>0$ 
sufficiently small. By virtue of equation (\ref{aux}), we
have
\begin{equation}
\label{aux1}
L_{0}u_{1}(x)=a(1-3w_{0}(x)){v_{1}^{2}(x)}-a{v_{1}^{3}(x)}-\varepsilon \sigma_{1}(x)
(w_{0}(x)+v_{1}(x)),
\end{equation}
\begin{equation}
\label{aux2}
L_{0}u_{2}(x)=a(1-3w_{0}(x)){v_{2}^{2}(x)}-a{v_{2}^{3}(x)}-\varepsilon \sigma_{1}(x)
(w_{0}(x)+v_{2}(x)).
\end{equation}
Formulas (\ref{aux1}) and (\ref{aux2}) yield
$$
L_{0}(u_{1}(x)-u_{2}(x))=(v_{1}(x)-v_{2}(x))\{a(1-3w_{0}(x))(v_{1}(x)+v_{2}(x))-
$$
\begin{equation}
\label{u12}
-a(v_{1}^{2}(x)+v_{1}(x)
v_{2}(x)+v_{2}^{2}(x))-\varepsilon \sigma_{1}(x)\}.
\end{equation}
We obtain the upper bound on the right side of equality (\ref{u12}) in the
absolute value as
$$
|v_{1}(x)-v_{2}(x)|\{a(1+3\|w_{0}\|_{L^{\infty}({\mathbb R}^{d})})(\|v_{1}\|_
{L^{\infty}({\mathbb R}^{d})}+\|v_{2}\|_{L^{\infty}({\mathbb R}^{d})})+
$$
\begin{equation}
\label{ub6}
a(\|v_{1}\|_{L^{\infty}({\mathbb R}^{d})}^{2}+
\|v_{1}\|_{L^{\infty}({\mathbb R}^{d})}\|v_{2}\|_{L^{\infty}({\mathbb R}^{d})}+
\|v_{2}\|_{L^{\infty}({\mathbb R}^{d})}^{2})+\varepsilon \|\sigma_{1}\|_
{L^{\infty}({\mathbb R}^{d})}\}.
\end{equation}
According to the Sobolev embedding (\ref{em}), quantity (\ref{ub6}) can be
estimated from above by
$$
|v_{1}(x)-v_{2}(x)|\{ac_{e}(1+3c_{e}\|w_{0}\|_{H^{2}({\mathbb R}^{d})})(\|v_{1}\|_
{H^{2}({\mathbb R}^{d})}+\|v_{2}\|_{H^{2}({\mathbb R}^{d})})+
$$
\begin{equation}
\label{ub7}
ac_{e}^{2}(\|v_{1}\|_{H^{2}({\mathbb R}^{d})}^{2}+
\|v_{1}\|_{H^{2}({\mathbb R}^{d})}\|v_{2}\|_{H^{2}({\mathbb R}^{d})}+
\|v_{2}\|_{H^{2}({\mathbb R}^{d})}^{2})+\varepsilon \|\sigma_{1}\|_
{L^{\infty}({\mathbb R}^{d})}\}.
\end{equation}
Let us recall that $v_1$, $v_2 \in B_{\rho}$. This yields the upper bound for
(\ref{ub7}) as
\begin{equation}
\label{ub8}
|v_{1}(x)-v_{2}(x)|\{2ac_{e}(1+3c_{e}\|w_{0}\|_{H^{2}({\mathbb R}^{d})})\rho+
3ac_{e}^{2}{\rho}^{2}+\varepsilon \|\sigma_{1}\|_{L^{\infty}({\mathbb R}^{d})}\}.
\end{equation}
Therefore, by means of (\ref{u12}) we arrive at
$$
\|u_{1}-u_{2}\|_ {H^{2}({\mathbb R}^{d})} \leq \|L_{0}^{-1}\|\{2ac_{e}(1+3c_{e}\|w_{0}\|
_{H^{2}({\mathbb R}^{d})})\rho+
$$
\begin{equation}
\label{u12n}
\begin{aligned}
+3ac_{e}^{2}{\rho}^{2}+\varepsilon \|\sigma_{1}\|_
{L^{\infty}({\mathbb R}^{d})}\}\|v_{1}-v_{2}\|_ {H^{2}({\mathbb R}^{d})}.
\end {aligned}
\end{equation}
Obviously, the estimate
\begin{equation}
\label{ub9}
\|L_{0}^{-1}\|\{2ac_{e}(1+3c_{e}\|w_{0}\|
_{H^{2}({\mathbb R}^{d})})\rho+3ac_{e}^{2}{\rho}^{2}+\varepsilon \|\sigma_{1}\|_
{L^{\infty}({\mathbb R}^{d})}\}<1
\end{equation} 
can be achieved for all $\rho>0$ and $\varepsilon>0$ small enough.
Hence, the map $t: B_{\rho}\to B_{\rho}$ defined by problem (\ref{aux}) is
a strict contraction. Its unique fixed point $w_{p}(x)$ is the only solution
of equation (\ref{pert}) in the ball $B_{\rho}$. Note that under the stated
assumptions the function $w_{p}(x)$ does not vanish identically for
$\varepsilon>0$. This can be easily seen from (\ref{pert}).
Evidently, the resulting solution
$w(x)$ of (\ref{stp}) given by formula (\ref{0p}) is contained in
$H^{2}({\mathbb R}^{d})$. Suppose the radius of the ball $B_{\rho}$ is
sufficiently small, such that
\begin{equation}
\label{ro}
\rho<\|w_{0}\|_{H^{2}({\mathbb R}^{d})}.
\end{equation}
By virtue of (\ref{0p}) and (\ref{ro}), using the triangle inequality
we obtain
\begin{equation}
\label{w}
\|w\|_{H^{2}({\mathbb R}^{d})}\geq \|w_{0}\|_{H^{2}({\mathbb R}^{d})}-
\|w_{p}\|_{H^{2}({\mathbb R}^{d})}\geq
\|w_{0}\|_{H^{2}({\mathbb R}^{d})}-\rho>0.
\end{equation}
Therefore, $w(x)$ is nontrivial as well.

Let us complete the proof  by obtaining
the asymptotics for the function $w_{p}(x)$ to the leading order in the
parameter $\varepsilon$. By means of (\ref{pert}), we have
$$
w_{p}(x)={L_{0}}^{-1}[a(1-3w_{0}(x)){w_{p}^{2}(x)}-
$$
\begin{equation}
\label{pert-}
-aw_{p}^{3}(x)-\varepsilon
\sigma_{1}(x)w_{p}(x)]-\varepsilon {L_{0}}^{-1}[\sigma_{1}(x)w_{0}(x)].
\end{equation}
Clearly, the leading term in the small parameter $\varepsilon$ in
the right side of (\ref{pert-}) equals to
\begin{equation}
\label{l}
-\varepsilon {L_{0}}^{-1}[\sigma_{1}(x)w_{0}(x)].
\end{equation}
Obviously, (\ref{l}) can be bounded from above in the $H^{2}({\mathbb R}^{d})$
norm by
\begin{equation}
\label{ln}
\varepsilon \|L_{0}^{-1}\|\|\sigma_{1}(x)\|_{L^{\infty}({\mathbb R}^{d})}\|w_{0}(x)\|_
{H^{2}({\mathbb R}^{d})}<\infty
\end{equation}
under the given conditions and via (\ref{l0-}). Evidently, the estimate from
above
$$
|a(1-3w_{0}(x)){w_{p}^{2}(x)}-aw_{p}^{3}(x)-\varepsilon \sigma_{1}(x)w_{p}(x)|\leq
a(1+3\|w_{0}\|_{L^{\infty}({\mathbb R}^{d})})\times
$$
\begin{equation}
\label{rub}
\times \|w_{p}\|_{L^{\infty}({\mathbb R}^{d})}|w_{p}(x)|+
a\|w_{p}\|_{L^{\infty}({\mathbb R}^{d})}^{2}|w_{p}(x)|+
\varepsilon \|\sigma_{1}\|_{L^{\infty}({\mathbb R}^{d})}|w_{p}(x)|
\end{equation}
is valid.
According to the Sobolev embedding (\ref{em}), the right side of 
(\ref{rub}) can be bounded from above by
$$
[ac_{e}(1+3c_{e}\|w_{0}\|_{H^{2}({\mathbb R}^{d})})\|w_{p}\|_{H^{2}({\mathbb R}^{d})}+
$$  
\begin{equation}
\label{rub1}
+ac_{e}^{2}\|w_{p}\|_{H^{2}({\mathbb R}^{d})}^{2}+\varepsilon \|\sigma_{1}\|_
{L^{\infty}({\mathbb R}^{d})}]|w_{p}(x)|.
\end{equation}
This means that the remaining term in the right side of (\ref{pert-}) can be
estimated from above in the $H^{2}({\mathbb R}^{d})$ norm by
$$
\|L_{0}^{-1}\|[ac_{e}(1+3c_{e}\|w_{0}\|_{H^{2}({\mathbb R}^{d})})
\|w_{p}\|_{H^{2}({\mathbb R}^{d})}+
$$
\begin{equation}
\label{rem}
+ac_{e}^{2}\|w_{p}\|_{H^{2}({\mathbb R}^{d})}^{2}+\varepsilon \|\sigma_{1}\|_
{L^{\infty}({\mathbb R}^{d})}]\|w_{p}\|_{H^{2}({\mathbb R}^{d})}=O({\varepsilon}^{2}).
\end{equation}
Therefore, identity (\ref{ass}) holds.  \hspace{6.5cm} $\Box$

\bigskip

%%%%%%%%%%%%%%%%%%%%%%%%%%%%%%%%%%%%%%%%%%%%%

\section{Acknowledgement} Stimulating discussions with Vitaly Volpert are
gratefully acknowledged.
The work was partially supported by the NSERC Discovery grant.

\bigskip

%%%%%%%%%%%%%%%%%%%%%%%%%%%%%%%%%%


\begin{thebibliography}{99}

\bibitem{AM07} A. ~Ambrosetti, A. Malchiodi, {\it Nonlinear analysis and
semilinear elliptic problems}, Cambridge Studies in Advanced Mathematics
{\bf 104}, Cambridge University Press, Cambridge (2007), 316 pp.  


\bibitem{BLP81} H. ~Berestycki, P.-L. ~Lions, L.A. ~Peletier, {\it An ODE
approach to the existence of positive solutions for semilinear problems in
${\mathbb R}^{N}$}, Indiana Univ. Math. J. {\bf 30} (1981), no. 1,  141--157.


\bibitem{CV21} Y. ~Chen, V. ~Vougalter, {\it Persistence of pulses for some
reaction-diffusion equations}, Pure Appl. Funct. Anal. {\bf 6} (2021), no. 2,
309--315.  

\bibitem{EV24} M. ~Efendiev, V. ~Vougalter,  {\it Solvability of some systems
of non-Fredholm integro-differential equations with mixed diffusion},
J. Dynam. Differential Equations {\bf 36} (2024), no. 3, 2239--2257.


\bibitem{EVV17} N. ~Eymard, V. ~Volpert, V. ~Vougalter, {\it Existence of pulses
for local and nonlocal reaction-diffusion equations}, J. Dynam. Differential
Equations {\bf 29} (2017), no. 3, 1145--1158.


\bibitem{IW00} D. ~Iron, M.J. ~Ward, {\it A metastable spike solution for a
nonlocal reaction-diffusion model}, SIAM J. Appl. Math. {\bf 60} (2000), no. 3,
778--802.


\bibitem{KP97}  I. ~Kuzin, S. ~Pohozaev, {\it Entire solutions of semilinear
elliptic equations}, Progress in Nonlinear Differential Equations and their
Applications {\bf 33}.  Birkh\"auser, Basel (1997), 250 pp.

  
\bibitem{MV23} M. ~Marion, V. ~Volpert, {\it Existence of pulses for
monotone reaction-diffusion systems}, SIAM J. Math. Anal. {\bf 55} (2023),
no. 2, 603--627.  


\bibitem{TBMV11}  J.C. ~Tzou, A. ~Bayliss,  B.J. ~Matkowsky,  V.A. ~Volpert,
{\it Stationary and slowly moving localised pulses in a singularly perturbed
Brusselator model}, European J. Appl. Math. {\bf 22} (2011), no. 5, 423--453.


\bibitem{V14} V. ~Volpert, {\it Elliptic partial differential equations},
Vol. 2. Reaction-diffusion equations. Monographs in Mathematics {\bf 104}.
Birkh\"auser/Springer, Basel (2014), 784 pp.


\bibitem{VV141} V. ~Volpert, V. ~Vougalter, {\it Existence of stationary pulses
for nonlocal reaction-diffusion equations}, Doc. Math. {\bf 19} (2014),
1141--1153.


\bibitem{V10} V. ~Vougalter, {\it On threshold eigenvalues and resonances
for the linearized NLS equation},  Math. Model. Nat. Phenom.  {\bf 5} (2010),
no. 4, 448--469.

\bibitem{V24} V. ~Vougalter, {\it Existence of solutions of some quadratic 
integral equations in dimensions two and three}, Pure Appl. Funct. Anal.
{\bf 9} (2024), no. 5, 1387--1396.

\bibitem{VV24} V. ~Vougalter, V. ~Volpert, {\it Solvability of some systems
of integro-differential equations in population dynamics depending on
the natality and mortality rates}, Arnold Math. J. {\bf 10} (2024), no. 1,
1--22.
  

\bibitem{WW07} J. ~Wei, M. ~Winter, {\it Existence, classification and stability
analysis of multiple-peaked solutions for the Gierer-Meinhardt system in
${\mathbb R}^{1}$}, Methods Appl. Anal. {\bf 14} (2007), no. 2, 119--163.


\bibitem{Z04} F. ~Zhang,  {\it Coexistence of a pulse and multiple spikes and
transition layers in the standing waves of a reaction-diffusion system},
J. Differential Equations {\bf 205} (2004), no. 1, 77--155.

\end{thebibliography}
\end{document}